\documentclass[a4paper,12pt]{article}
\newcounter{minutes}
\divide\time by 60
\newcounter{hours}
\multiply\time by 60
\addtocounter{minutes}{-\time}

\usepackage[english]{babel}
\usepackage[T1]{fontenc}
\usepackage{amsthm,enumerate,amsmath,amscd,amssymb}
\usepackage[dvips]{epsfig}
\usepackage[dvips]{graphicx}
\usepackage{icomma}

\newcommand{\arcsinh}{\,\textnormal{arcsinh}\,}
\newcommand{\arctanh}{\,\textnormal{arctanh}\,}
\newcommand{\comment}[1]{}

\swapnumbers

\theoremstyle{plain}

\newtheorem{theorem}[equation]{Theorem}

\newtheorem{lemma}[equation]{Lemma}

\theoremstyle{remark}

\newtheorem{remark}[equation]{Remark}

\numberwithin{equation}{section}

\begin{document}

\begin{center}
{\Large \bf On Jordan type inequalities for hyperbolic functions}
\end{center}

\begin{center}
{\large \bf R. Klén, M. Visuri and M. Vuorinen}
\end{center}

Department of Mathematics, University of Turku, FI-20014, Finland

\bigskip

\noindent {\bf Abstract.} This paper deals with some inequalities
for trigonometric and hyperbolic functions such as the Jordan
inequality and its generalizations. In particular, lower and upper
bounds for functions such as $(\sin x)/x$ and $x/\sinh x$ are
proved.

\medskip




\section{Introduction}

During the past several years there has been a great deal of
interest in trigonometric inequalities \cite{b1,b4,dz,jy,l,wsd,zs}. The
classical Jordan inequality \cite[p. 31]{m}
\begin{equation}\label{jordansinequality}
  \frac{2}{\pi}x \le \sin x \le x, \quad 0<x<\pi/2,
\end{equation}
has been in the focus of these studies and many refinements have
been proved for it by S.-H. Wu and H.M. Srivastava \cite{ws1,ws2},
X. Zhang, G. Wang and Y. Chu \cite{zwc}, J.-L. Li and Y.-L. Li
\cite{ll}, \cite{l}, S.-H. Wu and L. Debnath \cite{wd1,wd2,wd3},
A.Y. \"Ozban \cite{o}, F. Qi, D.-W. Niu and B.-N. Guo\cite{qng}, L. Zhu \cite{z1,z2,z3,z4,z5,z6,z7,z8,z9,z10,z11}, J. S\'andor
\cite{s1,s2}, \'A. Baricz and S. Wu \cite{bw1,bw2}, E. Neuman and J. S\'andor \cite{ns}, R. P. Agarwal, 
Y.-H. Kim and S. K. Sen \cite{aks}, D.-W. Niu, Z.-H. Huo, J. Cao and F. Qi \cite{nhcq}, 
W. Pan and L. Zhu \cite{pz} and F. Qi and B.-N. Guo \cite{qg}. For a long list of recent papers on this topic see
\cite{zs} and for an extensive survey see \cite{qng}. The proofs are based on familiar methods of calculus. In
particular, a method based on a l'Hospital type criterion for
monotonicity of the quotient of two functions from G. D. Anderson,
M. K. Vamanamurthy and M. Vuorinen \cite{avv1} is a key tool in
these studies. Some other applications of this criterion are
reviewed in \cite{avv2,avv3}. I. Pinelis has found several
applications of this criterion in \cite{p} and in several other
papers.

The inequality
\begin{equation}\label{baricz}
  \frac{1+\cos x}{2} \le \frac{\sin x}{x} \le \frac{2+\cos x}{3},
\end{equation}
where $x \in (-\pi/2,\pi/2)$ is well-known and it was studied
recently by \'A. Baricz in \cite[p. 111]{b2}. The second inequality of
(\ref{baricz}) is given in \cite[p. 354, 3.9.32]{m} for $0 \le x \le
\pi\,.$ For a refinement of the first inequality in (\ref{baricz})
see Remark \ref{refinement} (1) and of the second inequality see
Theorem \ref{refinement2}.

This paper is motivated by these studies and it is based on the Master Thesis of M. Visuri \cite{le}.
 Some of our main results are the following theorems.

\begin{theorem}\label{theorem1}
For $x \in \left( 0,\pi/2 \right)$
$$
  \frac{x^2}{\sinh ^2 x} < \frac{\sin x}{x} < \frac{x}{\sinh x}.
$$
\end{theorem}

\begin{theorem}\label{theorem2}
For $x \in (0,1)$
  $$
    \left( \frac{1}{\cosh x} \right) ^{\frac{1}{2}} <
    \frac{x}{\sinh x} < \left( \frac{1}{\cosh x} \right) ^{\frac{1}{4}}.
  $$
\end{theorem}

We will consider quotients $\sin x/x$ and $x/\sinh x$ at origin as limiting values
$\lim_{x \to 0} \sin x/x = 1$ and $\lim_{x \to 0}x/\sinh x = 1$.

\begin{remark}\label{refinement}
  (1) Let
  \[
    g_1(x) = \frac{1+2 \cos x}{3}+\frac{x \sin x}{6} \quad \text{and} \quad g_2(x) = \frac{1+\cos x}{2}.
  \]
  Then $g_1(x)-g_2(x) > 0$ on $(0,\pi/2)$ because
  \[
    \frac{d}{dx}(g_1(x)-g_2(x)) = \frac{x \cos x}{6} > 0.
  \]
  In \cite[(27)]{qcx} it is proved that for $x \in (0, \pi/2)$
  \[
    \frac{\sin x}{x} \ge g_1(x).
  \]
  Hence $(1+2 \cos x)/3+(x \sin x)/6$ is a better lower bound for $(\sin x)/x$ than (\ref{baricz}) for $x \in (0,\pi/2)$.

  (2) Observe that
  \begin{equation}\label{remark2}
    \frac{2+\cos x}{3} = \frac{2+2 \cos^2 (x/2)-1}{3} \le \cos \frac{x}{2}
  \end{equation}
  which holds true as equality if and only if $\cos (x/2) = (3 \pm 1)/4$. In conclusion,
  (\ref{remark2}) holds for all $x \in (-2\pi/3,2\pi/3)$. Together with (\ref{baricz}) we now have
  \[
    \cos^2 \frac{x}{2} = \frac{1+\cos x}{2} > \cos x
  \]
  and by (\ref{remark2})
  \[
    \cos^2\frac{x}{2} < \frac{\sin x}{x} < \cos \frac{x}{2}.
  \]
\end{remark}

\section{Jordan's inequality}

In this section we will find upper and lower bounds for $(\sin x)/x$ by using
hyperbolic trigonometric functions.

\begin{theorem}\label{thm1}
For $x \in (0,\pi/2)$
$$
  \frac{1}{\cosh x} < \frac{\sin x}{ x} < \frac{x}{\sinh x}.
$$
\end{theorem}
\begin{proof}
  The lower bound of $\sin x/x$ holds true if the function
  $f(x) = \sin x \cosh x -x$ is positive on $(0,\pi/2)$. Since
  $$
    f''(x) = 2\cos x \sinh x
  $$
  we have $f''(x) > 0$ for $x \in (0,\pi/2)$ and $f'(x)$ is increasing on $(0,\pi/2)$. Therefore
  $$
    f'(x) = \cos x \cosh x + \sin x \sinh x -1 > f'(0) = 0
  $$
  and the function $f(x)$ is increasing on $(0,\pi/2)$. Now $f(x) > f(0) = 0$ for $x \in (0,\pi/2)$.

  The upper bound of $\sin x/x$ holds true if the function $g(x) = x^2-\sin x \sinh x$ is positive on $(0,\pi/2)$. Let us denote $h(x) = \tan x -\tanh x$. Since $\cos x < 1 < \cosh x$ for $x \in (0,\pi/2)$ we have $h'(x) = \cosh^{-2}x-\cos^{-2}x > 0$ and $h(x) > h(0) = 0$ for $x \in (0,\pi/2)$. Now
  $$
    g'''(x) = 2(\cos x \cosh x)h(x),
  $$
  which is positive on $(0,\pi/2)$, because $\cos x \cosh x > 0$ and $h(x) > 0$ for $x \in (0,\pi/2)$. Therefore
  $$
    g''(x) = 2(1-\cos x \cosh x) > g(0) = 0
  $$
  and
  $$
    g'(x) = 2x- \cos x \sinh x - \sin x \cosh x > g'(0) = 0
  $$
  for $x \in (0,\pi/2)$. Now $g(x) > g(0) = 0$ for $x \in (0,\pi/2)$.
\end{proof}

\begin{proof}[Proof of Theorem \ref{theorem1}]
  The upper bound of $\sin x/x$ is clear by Theorem \ref{thm1}. The lower bound of $\sin x/x$ holds true if the function $f(x) = \sin x \sinh ^2 x-x^3$ is positive on $(0,\pi/2)$.

  Let us assume $x \in (0,\pi/2)$. Since $\sin x > x- x^3/6=(6x-x^3)/6$ we have $f(x) > ((6x-x^3) \sinh^2 x)/6-x^3$. We will show that
  \begin{equation}\label{gx}
    g(x) = \frac{6-x^2}{6} \sinh^2 x-x^2
  \end{equation}
  is positive which implies the assertion.

  Now $g(x) > 0$ is equivalent to
  $$
    \frac{\sinh x}{x} > \frac{\sqrt{6}}{\sqrt{6-x^2}}.
  $$
  Since $x^{-1} \sinh x > 1+x^2/6$ it is sufficient to show that $1+x^2/6 > \sqrt{6}/{\sqrt{6-x^2}}$, which is equivalent to
  \begin{equation}\label{poly}
    x^2(-x^4-6x^2+36) > 0.
  \end{equation}
  Let us denote $h(x) = -x^4-6x^2+36$. Now $h'(x) = -4x(x^2+3)$ and therefore $h'(x) \ne 0$ and $h(x) > \min \{ h(0),h( \pi/2) \} > 0$. Therefore inequality (\ref{poly}) holds for $x \in (0,\pi/2)$ and the assertion follows.
\end{proof}

We next show that for $x \in (0,1)$ the upper and lower bounds of (\ref{baricz}) are better than the upper and lower bounds in Theorem \ref{thm1}.

\begin{theorem}
  \renewcommand{\labelenumi}{\roman{enumi})}
  \begin{enumerate}
  \item For $x \in (-1,1)$
  \[
    \frac{2+\cos x}{3} \le \frac{x}{\sinh x}.
  \]
  \item For $x \in (-\pi/2,\pi/2)$
  \[
    \frac{1}{\cosh x} \le \frac{1+\cos x}{2} = \cos^2 \frac{x}{2}.
  \]
  \item For $x \in (-\pi/2,\pi/2)$
  \[
    \frac{1}{1+\sin^2 x} \le \frac{1+\cos^2 x}{2} \le \frac{1+\cos x}{2}.
  \]
  \end{enumerate}
\end{theorem}
\begin{proof}
  \textit{i)} The claim holds true if the function $f(x) = 3x-2\sinh x -\sinh x \cos x$ is non-negative on $[0,1)$. By a simple computation we obtain $f''(x) = 2(\cosh x \sin x -\sinh x)$. Inequality $f''(x) \ge 0$ is equivalent to $\sin x \ge \tanh x$. By the series expansions of $\sin x$ and $\tanh x$ we obtain
  \begin{eqnarray*}
    \sin x-\tanh x & = & \sum_{n \ge 3, \, n \equiv 1 (\textnormal{mod }2)} \frac{(-1)^{(n-1)/2}(n+1)-2^{n+1}(2^{n+1}-1)B_{n+1}}{(n+1)!}x^n\\
    & = & \sum_{n \ge 3, \, n \equiv 1 (\textnormal{mod }2)} c_n x^n,
  \end{eqnarray*}
  where $B_j$ is the $j\,$th Bernoulli number. By the properties of the Bernoullin numbers $c_1 = 1/6$, $c_3 = -1/8$, coefficients $c_n$, for $n \equiv 1 (2)$, form an alternating sequence, $|c_n x^n| \to 0$ as $n \to \infty$ and $|c_{2m+1}| > |c_{2m+3}|$ for $m \ge 1$. Therefore by Leibniz Criterion
  $$
    \sin x-\tanh x \ge \frac{x^3}{6}-\frac{x^5}{8} = \frac{x^3}{24} \left( 4-3x^2 \right)
  $$
  and $\sin x \ge \tanh x$ for all $x \in [0,1)$. Now $f(x)$ is a convex function on $[0,1)$ and $f'(x)$ is nondecreasing on $[0,1)$ with $f'(0)=0$. Therefore $f(x)$ is nondecreasing and $f(x) \ge f(0) = 0$.

  \textit{ii)} The claim holds true if the function $g(x) = \cosh x (1+\cos x)-2$ is non-negative on $[0,\pi/2)$. By the series expansion of $\cos x$ we have $\cos x-1+x^2/2 \ge 0$ and therefore by the series expansion of $\cosh x$
  \begin{eqnarray*}
    g(x) & \ge & \left( 1+\frac{x^2}{2} \right) (1+\cos x)-2\\
    & = & \cos x-1+\frac{x^2}{2}+\frac{x^2 \cos x}{2}\\
    & \ge & \cos x-1+\frac{x^2}{2} \ge 0
  \end{eqnarray*}
  and the assertion follows.

  \textit{iii)} Clearly we have
  \[
    (1+\cos^2 x)(1+\sin^2 x) = 2+\sin^2 x \cos^2 x \ge 2,
  \]
  which implies the first inequality of the claim. The second inequality is trivial since $\cos x \in (0,1)$.
\end{proof}

\begin{theorem}
  Let $x \in (0,\pi/2)$. Then
  \renewcommand{\labelenumi}{\roman{enumi})}
  \begin{enumerate}
  \item the function
  \[
    f(t) = \cos^t\frac{x}{t}
  \]
  is increasing on $(1,\infty)$,
  \item the function
  \[
    g(t) = \sin^t\frac{x}{t}
  \]
  is decreasing on $(1,\infty)$,
  \item the functions $\overline{f}(t) = \cosh^t(x/t)$ and $\overline{g}(t) = \sinh^t(x/t)$ are decreasing on $(0,\infty)$.
  \end{enumerate}
\end{theorem}
\begin{proof}
  \textit{i)} Let us consider instead of $f(x)$ the function
  \[
    f_1(y) = \frac{x}{y}\log \cos y,
  \]
  for $y \in (0,x)$. Note that $f(t) = \exp(f_1(x/t))$ and therefore the claim is equivalent to the function $f_1(y)$ being decreasing on $(0,x)$. We have
  \[
    f_1'(y) = -\frac{x}{y^2}(\log \cos y+y \tan y)
  \]
  and $f_1'(y) \le 0$ is equivalent to $f_2(y) = \log \cos y+y \tan y \ge 0$. Since $f_2'(y) = y/\cos^2 y \ge 0$ we have $f_2(y) \ge f_2(0) = 0$. Therefore $f(t)$ is increasing on $(1,\infty)$.

  \textit{ii)} We will consider instead of $g(x)$ the function
  \[
    g_1(y) = \frac{x}{y}\log \sin y,
  \]
  for $y \in (0,x)$. Note that $g(t) = \exp(g_1(x/t))$ and therefore the claim is equivalent to the function $g_1(y)$ being increasing on $(0,x)$. We have
  \[
    g_1'(y) = \frac{x}{y^2}(y / \tan y-\log \cos y)
  \]
  and $g_1'(y) \ge 0$ is equivalent to $g_2(y) = y / \tan y-\log \cos y \ge 0$. Since $g_2'(y) = ((1/\cos y) - (y/\sin y))/\sin y \ge 0$ we have $g_2(y) \ge f_2(0) = 1$. Therefore $g_1'(y) \ge 0$ and the assertion follows.

  \textit{iii)} We will show that $h_1(y) = (x/y) \log \cosh y$ is increasing on $(0,\infty)$. Now $h_1'(y) = (x (y \tanh y-\log \cosh y))/y^2$,
  \[
    \frac{d (y \tanh y-\log \cosh y)}{d y} = \frac{y}{\cosh^2 y} > 0
  \]
  and $y \tanh y-\log \cosh y \ge 0$. Therefore the function $h_1(y)$ is increasing on $(0,\infty)$ and $\overline{f}(t)$ is decreasing on $(0,\infty)$.

  We will show that $h_2(y) = (x/y) \log \sinh y$ is increasing on $(0,\infty)$. Now $h_2'(y) = (x (y \coth y-\log \sinh y))/y^2$,
  \[
    \frac{d (y \coth y-\log \sinh y)}{d y} = -\frac{y}{\sinh^2 y} < 0
  \]
  and $\coth y-(\log \sinh y)/y \ge \lim_{y \to \infty} \coth y-(\log \sinh y)/y = 0$. Therefore the function $h_2(y)$ is increasing on $(0,\infty)$ and $\overline{g}(t)$ is decreasing on $(0,\infty)$.
\end{proof}

We next will improve the upper bound of (\ref{baricz}).
\begin{theorem}\label{refinement2}
  For $x \in (-\sqrt{27/5},\sqrt{27/5})$
  \begin{equation}\label{cos3}
    \cos^2 \frac{x}{2} \le \frac{\sin x}{x} \le \cos^3 \left( \frac{x}{3} \right) \le \frac{2+\cos x}{3}.
  \end{equation}
\end{theorem}
\begin{proof}
  The first inequality of (\ref{cos3}) follows from (\ref{baricz}).

  By the series expansions of $\sin x$ and $\cos x$
  \[
    \frac{\sin x}{x} \le 1- \frac{x^2}{6}+\frac{x^4}{120} \le \left( 1-\frac{x^2}{18} \right)^3 \le \cos^3 \left( \frac{x}{3} \right),
  \]
  where the second inequality is equivalent to $x^4(27-5x^2)/29160 \ge 0$ and the second inequality of (\ref{cos3}) follows.

  By the identity $\cos^3 x = (\cos 3x+3\cos x)/4$ the upper bound of (\ref{cos3}) is equivalent to $0 \le 8+ \cos x - 9 \cos (x/3)$. By the series expansion of $\cos x$
  \[
    8+ \cos x - 9 \cos (x/3) = \sum_{n=2}^\infty (-1)^n\frac{3^{2n}-9}{3^{2n}(2n)!}x^{2n}
  \]
  and by the Leibniz Criterion the assertion follows.
\end{proof}

\section{Hyperbolic Jordan's inequality}

In this section we will find upper and lower bounds for the functions $x/\sinh x$ and $\cosh x$.

\begin{theorem}\label{hjordan}
  For $x \in (-\pi/2,\pi/2)$
  \[
      1-\frac{x^2}{6} \le \frac{\sin x}{x} \le 1- \frac{2x^2}{3\pi^2}.
  \]
\end{theorem}
\begin{proof}
  We obtain from the series expansion of $\sin x$
  \[
    \frac{\sin x}{x} \ge 1-\frac{x^2}{6},
  \]
  which proves the lower bound.

  By using the identity $1-\cos x = 2\sin^2 (x/2)$ the chain of inequalities (\ref{baricz}) gives
  \[
    \frac{\sin x}{x} \le 1- \frac{2\sin^2 (x/2)}{3}
  \]
  and the assertion follows from inequality $\sin^2 (x/2) \ge (x/\pi)^2$.
\end{proof}

\begin{remark}
  Li -Li have proved \cite[(4.9)]{ll} that
  $$
    \frac{\sin x}{x}  < p(x) \left( 1 -  \frac{x^2}{\pi^2} \right)  < 1-\frac{x^2}{\pi^2},\quad  0 < x < \pi,
  $$
  where  $p(x)= 1/\sqrt{1+ 3(x/ \pi)^4 }<1\,.$ This result improves Theorem \ref{hjordan}.
\end{remark}

\begin{lemma}\label{ubsinh}
  For $x \in (0,1)$
  \renewcommand{\labelenumi}{\roman{enumi})}
  \begin{enumerate}
    \item $\displaystyle \sinh x < x+\frac{x^3}{5},$
    \item $\displaystyle \cosh x < 1+\frac{5x^2}{9},$
    \item $\displaystyle \frac{1}{\cosh x} < 1-\frac{x^2}{3}.$
  \end{enumerate}
\end{lemma}
\begin{proof}
  \textit{i)} For $x \in (0,1)$ we have $x^2(1-x^2) > 0$ which is equivalent to
  \begin{equation}\label{simplealgebra}
    \frac{1}{1-x^2/6} < 1+\frac{x^2}{5}.
  \end{equation}
  By Theorem \ref{thm1}, Theorem \ref{hjordan} and (\ref{simplealgebra})
  $$
    \sinh x \le \frac{x^2}{\sin x}\le \frac{x}{1-x^2/6} < x+\frac{x^3}{5}.
  $$

  \textit{ii)} Since $(2n)! > 6^n$ for $n \ge 3$ we have
  \begin{eqnarray*}
     1+\frac{5x^2}{9}-\cosh x & = & \frac{x^2}{18} - \sum_{n=2}^\infty \frac{x^{2n}}{(2n)!}\\
     & \ge & \frac{x^2}{18} - \sum_{n=2}^\infty \frac{x^2}{(2n)!}\\
     & = & x^2 \left( \frac{1}{72} - \sum_{n=3}^\infty \frac{1}{(2n)!} \right)\\
     & \ge & x^2 \left( \frac{1}{72} - \sum_{n=3}^\infty \frac{1}{6^n} \right)\\
     & > & 0.
  \end{eqnarray*}

  \textit{iii)} By the series expansion of $\cosh x$ we have
  $$
    \cosh x \left( 1-\frac{x^2}{3} \right) \ge \left( 1+\frac{x^2}{2} \right) \left( 1-\frac{x^2}{3} \right) = 1+\frac{x^2}{6}-\frac{x^4}{6} > 1.
  $$
\end{proof}

\begin{proof}[Proof of Theorem \ref{theorem2}]
  The lower bound of $x/\sinh x$ follows from Lemma \ref{ubsinh} and Theorem \ref{hjordan} since
  $$
    \frac{1}{\cosh x} < 1-\frac{x^2}{3} < \left( 1-\frac{x^2}{6} \right)^2 \le \left( \frac{x}{\sinh x} \right)^2.
  $$

  The upper bound of $x/\sinh x$ holds true if the function $g(x) = \sinh^4 x -x^4 \cosh x$ is positive on $(0,1)$. By the series expansion it is clear that
  \begin{equation}\label{lbsinh}
    \sinh x > x+x^3/6.
  \end{equation}
  By Lemma \ref{ubsinh} and (\ref{lbsinh})
  $$
    g(x) > \left( x+\frac{x^3}{6} \right)^4 - x^4 \left( 1+\frac{2x^2}{3} \right) = \frac{x^6}{1296}(x^6+24x^4+216x^2+144) > 0
  $$
  and the assertion follows.
\end{proof}

\begin{theorem}\label{coshbound1}
  For $x \in (0,\frac{\pi}{4})$
  $$
    \cosh x < \frac{\cos x}{\sqrt{(\cos x)^2-(\sin x)^2}}.
  $$
\end{theorem}
\begin{proof}
  The upper bound of $\cosh x$ holds true if the function $f(x) = \cos^2 x-\cosh^2 x (\cos^2 x-\sin^2 x)$ is positive on $(0,\pi/4)$. Since
  $$
    f''(x) = 4 \sin (2x) \sinh (2x) >0
  $$
  we have
  $$
    f'(x) = \sin (2x) \sinh (2x)-\cos (2x) \cosh (2x) > f'(0) = 0.
  $$
  Therefore $f(x) > f(0) = 0$ and the assertion follows.
\end{proof}

\begin{theorem}\label{coshbound2}
  For $x \in (0,\frac{\pi}{4})$
  $$
    \frac{1}{(\cos x)^{2/3}} < \cosh x < \frac{1}{\cos x}.
  $$
\end{theorem}
\begin{proof}
  The upper bound of $\cosh x$ holds true if the function $f(x) = 1 - \cos x \cosh x$ is positive on $(0,\pi/4)$. Since $f''(x) = 2\sin x \sinh x > 0$ the function $f'(x) = \cosh x\sin x -\cos x \sinh x$ is increasing. Therefore $f'(x) > f'(0) = 0$ and $f(x) > f(0) = 0$.

  The lower bound of $\cosh x$ holds true if the function $g(x) = \cos^2 x \cosh^3 x - 1$ is positive on $(0,\pi/4)$. By the series expansions we have
  $$
    g(x) > \left( 1-\frac{x^2}{2} \right)^2 \left( 1+\frac{x^2}{2} \right)^3 -1 = \frac{x^2}{32}(x^8+2x^6-8x^4-16x^2+16).
  $$
  By a straightforward computation we see that the polynomial $h(x) = x^8+2x^6-8x^4-16x^2+16$ is strictly decreasing on $(0,\pi/4)$. Therefore
  \begin{eqnarray*}
    h(x) & > & h(\pi/4)\\
    & = & 16-\pi^2-\frac{\pi^4}{32}+\frac{\pi^6}{2048}+\frac{\pi^8}{65536}\\
    & > & 16-\frac{16}{5}^2-32^{-1}\frac{16}{5}^4+\frac{3^6}{2048}+\frac{3^8}{65536}\\
    & = & \frac{120392497}{40960000} > 0
  \end{eqnarray*}
  and the assertion follows.
\end{proof}

\begin{remark}
  \'A. Baricz and S. Wu have shown in \cite[pp. 276-277]{bw2} that the right hand side of Theorem \ref{thm1} is true for $x \in (0,\pi)$ and the right hand side of Theorem \ref{coshbound2} is true for $x \in (0,\pi/2)$. Their proof is based on the infinite product representations.
\end{remark}

Note that for $x \in (0,\pi/4)$
\[
  \frac{1}{\cos x} \le \frac{\cos x}{\sqrt{(\cos x)^2-(\sin x)^2}}.
\]
Hence, the upper bound in Theorem \ref{coshbound2} is better that in Theorem \ref{coshbound1}.

\section{Trigonometric inequalities}

\begin{theorem}\label{inverseineq}
  For $x \in (0,1)$ the following inequalities hold
  \renewcommand{\labelenumi}{\roman{enumi})}
  \begin{enumerate}
    \item $\displaystyle \frac{x}{\arcsin x} \le \frac{\sin x}{x},$
    \item $\displaystyle \frac{x}{\arcsinh x} \le \frac{\sinh x}{x},$
    \item $\displaystyle \frac{x}{\arctan x} \le \frac{\tan x}{x},$
    \item $\displaystyle \frac{x}{\arctanh x} \le \frac{\tanh x}{x}.$
  \end{enumerate}
\end{theorem}
\begin{proof}
  \textit{i)} By setting $x = \sin t$ the assertion is equivalent to
  \[
    \textrm{sinc} \, t \le \textrm{sinc} \, (\sin t),
  \]
  which is true because $\textrm{sinc}\,t = (\sin t)/t$ is decreasing on $(0,\pi/2)$ and $\sin t \le t$.


  \textit{ii)} By the series expansions of $\sinh x$ and $\arcsinh x$ we have by Leibniz Criterion
  \begin{eqnarray*}
    (\sinh x) \arcsinh x - x^2 & \ge & \left( x+\frac{x^3}{6} \right) \left( x-\frac{x^3}{6}+\frac{3x^5}{40}-\frac{5x^7}{112} \right)-x^2\\
    & = & \frac{x^6}{10080}(-75x^4-324x^2+476)
  \end{eqnarray*}
  and since \( -75x^4-324x^2+476 >77 \) on \( (0,1) \) the assertion follows.

  \textit{iii)} By the series expansions of $\tan x$ and $\arctan x$ we have by Leibniz Criterion
  \begin{eqnarray*}
    (\tan x) \arctan x - x^2 & \ge & \left( x+\frac{x^3}{3}+\frac{2x^5}{15}+\frac{17x^7}{315} \right) \left( x-\frac{x^3}{3} \right)-x^2\\
    & = & \frac{x^6}{945}(21+9x^2-17x^4)
  \end{eqnarray*}
  and since \( 21+9x^2-17x^4 > 4 \) on \( (0,1) \) the assertion follows.

  \textit{iv)} By the series expansions of $\tanh x$ and $\arctanh x$ we have by Leibniz Criterion
  \begin{eqnarray*}
    (\tanh x) \arctanh x - x^2 & \ge & \left( x-\frac{x^3}{3} \right) \left( x+\frac{x^3}{3}+\frac{x^5}{5} \right)-x^2\\
    & = & \frac{x^6}{45}(4-3x^2)
  \end{eqnarray*}
  and since \( 4-3x^2 > 1 \) on \( (0,1) \) the assertion follows.
\end{proof}

\begin{remark}
  Similar inequalities to Theorem \ref{inverseineq} have been considered by E. Neuman in \cite[pp. 34-35]{n}.
\end{remark}

\begin{theorem}
  Let $k \in (0,1)$. Then
  \renewcommand{\labelenumi}{\roman{enumi})}
  \begin{enumerate}
    \item for $x \in (0,\pi)$
    $$
      \frac{\sin x}{x} \le \frac{\sin (kx)}{kx},
    $$
    \item for $x > 0$
    $$
      \frac{\sinh x}{x} \ge \frac{\sinh (kx)}{kx},
    $$
    \item for $x \in (0,1)$
    $$
      \frac{\tanh x}{x} \le \frac{\tanh (kx)}{kx}.
    $$
  \end{enumerate}
\end{theorem}
\begin{proof}
  \textit{i)} The claim follows from the fact that $\textrm{sinc}$ is decreasing on $(0,\pi)$.

  \textit{ii)} The claim is equivalent to saying that the function $f(x) = (\sinh x)/x$ is increasing for $x > 0$. Since $f'(x) = (\cosh x)/x-(\sinh x)/x^2 \ge 0$ and $f'(x) \ge 0$ is equivalent to $\tanh x \le x$ the assertion follows.

  \textit{iii)} The claim is equivalent to $\tanh (kx) - k \tanh x \ge 0$. By the series expansion of $\tanh x$ we have
  \[
    \tanh (kx) - k \tanh x = k \sum_{n=1}^\infty \frac{4^{n+1}(4^{n+1}-1)B_{2(n+1)}x^{2n+1}}{(2n+2)!}(k^{2n}-1),
  \]
  where $B_j$ is the $j$th Bernoulli number ($B_0 = 1$, $B_1 = -1/2$, $B_2=1/6, \dots$). The assertion follows from the Leibniz Criterion, if
  \begin{equation}\label{leibniz2}
    \frac{k-k^3}{3}x^3-\frac{2(k-k^5)}{15}x^5 > 0
  \end{equation}
  for all $x \in (0,1)$. Since (\ref{leibniz2}) is equivalent to
  \[
    x^2 < \frac{5}{2(1+k^2)}
  \]
  the assertion follows from the assumptions $k \in (0,1)$ and $x \in (0,1)$.
\end{proof}


\end{document}